\documentstyle[12pt]{article} 
\def\1{\hbox{ 1\kern-.35em\hbox{1}}}
\newcommand{\Bbb}{\bf}
\newcommand{\frak}{\bf} 
\newcommand{\Z}{{\Bbb Z}} 
\newcommand{\C}{{\Bbb C}} 
\newcommand{\E}{{\Bbb E}} 
\newcommand{\Nn}{{\Bbb N}_n} 
\newcommand{\cP}{{\cal P}} 
\newcommand{\cT}{{\cal T}_q({\frak g})} 
\newcommand{\Tq}{{\cal T}_q({\frak g})} 
\newcommand{\ca}{{\cal T}_q({\frak g})}
\newcommand{\Gamq}{\Gamma_q^{\frak k}} 
\newcommand{\Hq}{\Gamq} 
\newcommand{\Eq}{{\cal A}_q^{\frak k}}

\newcommand{\osp}{\mbox{osp}(1|2n)} 
\newcommand{\OSPq}{\mbox{OSP}_q(1|2n)} 

\newcommand{\g}{{\frak g}} 
\newcommand{\Uq}{\mbox{U}_q({\frak g})} 
 
\newcommand{\Uk}{\mbox{U}_q({\frak k})} 
\newcommand{\Up}{\mbox{U}_q({\frak p})} 
\newcommand{\id }{\mbox{id}}

\newcommand{\irrep}{ irreducible representation} 
 
\newcommand{\Mod}{ \mbox{\bf Mod}_q(\g)} 
\newcommand{\ba}{\begin{eqnarray}}
\newcommand{\na}{\end{eqnarray}}
\newcommand{\ban}{\begin{eqnarray*}}
\newcommand{\nan}{\end{eqnarray*}}
\newcommand{\bea}{\ba} 
\newcommand{\eea}{\na}
\newcommand{\be}{\ban}
\newcommand{\ee}{\nan}  
\newtheorem{lemma}{Lemma}
\newtheorem{proposition}{Proposition}
\newtheorem{theorem}{Theorem}

\newcommand{\co}{\mbox{$\cal O$}}

\newenvironment{proof}[1]{\begin{trivlist} \item[] {\em #1\/}: }%
{\hfill $\Box$ \end{trivlist}}

\newcommand{\SDq}{\mbox{SD}_q}

\begin{document}
\hfill{\small NCU/CCS-1998-0330}
\vspace{1cm}
\large
{\bf \begin{center}
GEOMETRY AND REPRESENTATIONS OF THE\\ 
QUANTUM SUPERGROUP $\mbox{OSP}_q(1|2n)$
\end{center} }

\normalsize 
\begin{center} 
H. C. Lee$^1$ \  \  and \  \   R. B. Zhang$^2$\\                     
\end{center}   

\small 
\begin{center} 
$^1$ Department of Physics and Center for Complex Systems,\\             
National Central University, Chungli, Taiwan, ROC. \\             
$^2$ Department of Pure Mathematics, \\             
University of Adelaide, Adelaide, Australia.
\end{center}

\normalsize 
\begin{abstract}
The quantum supergroup $\OSPq$ is studied  systematically. 
A Haar functional is constructed, and an 
algebraic version of the Peter - Weyl theory is 
extended to this quantum supergroup. 
Quantum homogeneous superspaces and quantum homogeneous supervector 
bundles are defined following the strategy of Connes' theory.  
Parabolic induction is developed 
by employing the quantum homogeneous supervector bundles. 
Quantum Frobenius reciprocity and a generalized Borel - Weil 
theorem are established for the induced representations. 
\end{abstract}
 
\vspace{1cm}

\section{\normalsize{\bf INTRODUCTION}}
Quantized universal enveloping algebras of Lie superalgebras 
were introduced in the late 1980s \cite{Bracken, Chaichian}
to describe the type of supersymmetries exhibited 
by some $2$ - dimensional statistical mechanics models \cite{Perk}.
Since then these quantum superalgebras have been intensively 
studied, leading to the development of an extensive theory 
on both the structure and representations. 
We mention in particular that the quasi - triangular Hopf superalgebraic
structure of the quantum superalgebras was investigated in \cite{Double};
the representation theory of the type I 
quantum superalgebras,  the $gl(m|n)$  super Yangians
and the quantum affine superalgebras with symmetrizable 
Cartan matrices were developed in 
\cite{II}. The theory of quantum superalgebras 
had significant impact on  
a range of areas of physics and mathematics. 
Its applications to two dimensional
integrable models in statistical mechanics and quantum
field theory were extensively explored in  
\cite{Bracken, Delius} and many other publications.
The  application to knot theory and $3$ - manifolds 
\cite{Links}\cite{Manifolds} has 
yielded many new topological invariants, notably, the
multi - parameter generalizations of Alexander - Conway polynomials.

The associated quantum supergroups are in contrast less studied 
in the literature.  So far only the quantum supergroup 
$GL_q(m|n)$ has been systematically investigated \cite{Zhang}.  
In \cite{Zhang}, the structure and representation theories    
of $GL_q(m|n)$ were developed.  The irreducible 
covariant and contravariant tensorial  representations were 
studied in detail within the framework of parabolic induction, 
resulting in a quantum Borel - Weil theorem for  
these representations. The aim of this paper is to  
treat the $\osp$  series of quantum supergroups at generic $q$. 

The $\osp$ - series of Lie superalgebras, $\mbox{osp}(1|32)$ in 
particular, featured prominently in recent developments 
of string theory.  An Inonu - Wigner contraction of 
$\mbox{osp}(1|32)$ yields the 
11 - dimensional Poincar\'e superalgebra with two and five form 
central charges, which is the underlying symmetry of M - theory; 
the superalgebra $\mbox{osp}(1|32)$ itself also plays an important 
role in the theory of supermembranes \cite{Townsend}. 
From a mathematical point of view, 
$\osp$ is also rather exceptional amongst all the finite 
dimensional simple Lie superalgebras 
in that its Cartan matrix is symmetrizable, and the structure 
of its finite dimensional representations   is completely understood. 
In particular, it is known that all finite dimensional 
representations are completely reducible. 

Many properties of $\osp$ carry over to the quantum case when $q$ 
is generic.  It is particularly useful to recall that the Drinfeld 
version of $\mbox{U}_q(\osp)$ is, algebraically, 
a trivial deformation of $\mbox{U}(\osp)$\footnote{This fact is 
known to experts, and may be easily inferred from results 
of \cite{Scheunert}.} 
in the sense of Gerstenhaber.  Therefore, 
{\em finite dimensional representations of} $\mbox{U}_q(\osp)$ 
{\em are also completely reducible}. 
This remains true for the Jimbo version of $\mbox{U}_q(\osp)$ 
at generic $q$. 
One way to see this is through the specialization of 
the indeterminate of the Drinfeld algebra to a generic 
complex parameter; the other is through the 
isomorphism between $\mbox{U}_q(\osp)$ and $\mbox{U}_{-q}(so(2n+1))$
established by a kind of Bose - Fermi 
transmutation \cite{I}.  There is also an interesting 
connection between the representation theory of 
$\mbox{U}_q(\osp)$  and quantum para - statistics, details on 
which can be found in \cite{Palev}.

This paper will study structural and representation theoretical 
properties of the quantum supergroup $\OSPq$,  and also investigate 
its underlying geometries.  This quantum supergroup will be defined 
by its superalgebra of functions, which is the $\Z_2$ - graded Hopf 
algebra generated by the matrix elements of the vector 
representation of $\mbox{U}_q(\osp)$.  Two major results in the structure 
theory are presented, namely, 
the existence of a left and right integral, 
which will be called a quantum Haar functional, 
and a quantum Peter - Weyl theorem.  
 
Corresponding to each reductive subalgebra $\Uk$  of $\mbox{U}_q(\osp)$, 
we introduce a quantum homogeneous superspace, which is 
defined by specifying its superalgebra of functions $\Eq$.  
A quantum homogeneous supervector bundle over the 
quantum homogeneous superspace is induced from any given 
finite dimensional $\Uk$ - module. We shall show that 
the space of sections $\Hq(V)$ of this bundle is projective 
and is of finite type both as a left and a right module 
over $\Eq$. Therefore our definition of 
quantum homogeneous supervector bundles is consistent
with the general definition of noncommutative vector bundles 
in Connes' theory \cite{Connes}.    

Quantum homogeneous supervector bundles will be applied to 
develop a theory of induced representations for $\OSPq$. 
Amongst the results obtained are quantum versions of 
Frobenius reciprocity and the Borel - Weil theorem. 
The latter provides a concrete realization of finite 
dimensional irreducible $\OSPq$ representations in terms 
of quantum analogs of `holomorphic' sections of quantum
homogeneous supervector bundles.  

We wish to point out that in the context of Lie supergroups 
at the classical level, the mathematical theories of 
homogeneous superspaces and 
homogeneous supervector bundles were studied in 
\cite{Manin, Penkov}.  The development of 
a Bott - Borel - Weil theory was also initiated and extensively
investigated by Penkov and co - workers\cite{Penkov}. 
However, complications arising from supermanifold  geometry 
render these subjects very difficult to study.  So far as 
we are aware, many aspects of the subjects remain to be fully 
developed.  It seems that the Hopf algebraic approach 
developed here and in \cite{Zhang} is also worth exploring
at the classical level, and is likely to
provide a new method complementary to the geometric approach
of \cite{Manin, Penkov}.

The organization of the paper is as follows. 
Section 2 reviews some known facts about $\mbox{U}_q(\osp)$, 
which will be needed later. Section 3 studies the quantum 
supergroup $\OSPq$.  Section 4 investigates the quantum 
homogeneous superspaces and quantum homogeneous 
supervector bundles determined by this quantum supergroup, 
while the last section applies results of section 4 to 
study the representation theory of $\OSPq$.

\section{\normalsize{$\mbox{U}_q(\mbox{osp}(1|2n))$}}
This section reviews some known results on the 
quantized universal enveloping algebra $\mbox{U}_q(\mbox{osp}(1|2n))$. 
Let $E$ be the $n$ - dimensional Euclidean space spanned by 
the vectors $\epsilon_i$, with the inner product $(\, ,\, )$ 
defined by $(\epsilon_i,\, \epsilon_j)=\delta_{i j}$. 
We can express the simple roots of the Lie superalgebra 
$\osp$ in terms of the $\epsilon$'s as 
\be 
\alpha_i& =&\epsilon_i -\epsilon_{i+1}, \quad i=1, 2, ..., n-1, \\
\alpha_n&=&\epsilon_n,  
\ee     
where $\alpha_n$ is the odd simple root. 
The Cartan matrix $A=\left(a_{i j}\right)_{i, j=1}^n$ 
of $\osp$ is then given by  
$a_{ i j}={{2(\alpha_i,\,  \alpha_j)}\over{(\alpha_i,\,  \alpha_i)} }$.
An element $\mu\in E$ will be called integral if   
\be 
l_i&=&{{2(\mu,\,  \alpha_i)}\over {(\alpha_i,\,  \alpha_i)} }\in\Z, 
\quad \forall i<n, \\
l_n&=&{{(\mu,\,  \alpha_n)}\over {(\alpha_n,\,  \alpha_n)} }\in\Z, 
\ee 
and the set of all integral elements will be denoted by $\cP$.
(Note the unusual form of $l_n$.)  
Set $\cP_+ =$  $\{\mu\in\cP\, |\, l_i, l_n\in\Z_+\}$. Elements of 
$\cP_+$ will be called integral dominant.

The Jimbo version of the quantum superalgebra 
$\mbox{U}_q(\osp)$ is a $\Z_2$ - graded complex associative 
algebra generated by $\{  k_i^{\pm 1}, \ e_i, \ f_i, 
\ i\in\Nn\}$, $\Nn=\{1, 2, ..., n\}$,
subject to the relations
\bea
k_i k^{-1}_i=1, &   k_i k_j = k_j k_i,  \nonumber \\
k_i e_j = q^{(\alpha_i, \ \alpha_j)} e_j k_i,
&k_i f_j= q^{-(\alpha_i, \ \alpha_j)} f_j k_i,\nonumber \\
{[}e_i, f_j\} = \delta_{i j}{{k_i - k_i}\over
{q - q }},
&\forall i, j\in I, \nonumber \\
(Ad e_i)^{1-a_{i j}}(e_j) = 0, &(Ad f_i)^{1-a_{i j}}(f_j)=0,
\ \   \forall i\ne j. 
\label{definition}
\eea
All the generators are chosen to be homogeneous,
with  $k^{\pm 1}_i, \ \forall i$, and $e_j$, $f_j$, $j<n$,
being even, and $e_n$, $f_n$ being odd.
For a homogeneous element $x$, we define $[x]=0$ if $x$ is even,
and $[x]=1$ when odd.  The graded commutator $[ .\  ,\ . \}$
represents the usual commutator when any one of the two arguments
is even, and the anti commutator when both arguments are odd.
The adjoint operation $Ad$ is defined by
\ban
Ad e_i(x)&=& e_i x - (-1)^{[e_i][x]} k_i x k_i^{-1} e_i, \\
Ad f_i(x)&=& f_i x - (-1)^{[f_i][x]} k_i x k_i^{-1} f_i.
\nan
For $x$ being a monomial in $e_j$'s or $f_j$'s, it carries a
definite weight $\omega(x)\in H^*$. Then
$Ad e_i(x)$ $=$ $e_i x $ $-$
$(-1)^{[e_i][x]} q^{(\alpha_i, \ \omega(x))} x e_i$,
and similarly for $Ad f_i(x)$.  For convenience, we will 
use the notation $\frak g$ to denote $\osp$, and 
$\Uq$ to denote $\mbox{U}_q(\osp)$.
As is well known, this algebra has the structures of a ${\bf Z}_2$ 
graded Hopf algebra. We will denote the co-multiplication by 
$\Delta$, the co-unit by $\epsilon$ and the antipode by $S$.

The representation theory of $\Uq$ was developed in \cite{I}.  
For any finite dimensional $\Uq$ - module, there exists 
a homogeneous basis relative to which  the $k_i$ are represented 
by diagonal matrices. Here we will only consider such finite 
dimensional $\Uq$ - modules
that the eigenvalues of the $k_i$ tend to $1$ as $q$ approaches $1$.
We will denote the set of all such $\Uq$ - modules by $\Mod$.  
Recall that all objects of $\Mod$ are semi-simple.

If $W(\lambda)$ is a simple object of $\Mod$,  
then there exists the unique (up to scalar multiples)
highest weight  vector $v_+$, such that
\ban
e_i v_+ =0, & \quad  k_i v_+ = q^{(\lambda, \ \alpha_i)} v_+,
& \quad \lambda\in\cP_+,
\nan
and the module $W(\lambda)$ is uniquely determined by
the highest weight $\lambda$.
We will denote the lowest weight of
$W(\lambda)$ by $\bar\lambda$, and define
$\lambda^\dagger = -\bar\lambda$.
The dual module of $W(\lambda)$ has highest weight $\lambda^\dagger$.

The irreducible $\Uq$ - module with highest weight $\epsilon_1$ 
plays a special role in the representation theory of $\Uq$. 
We denote this module by  $\E$, and refer to it as the vector module.
Let us now examine this module in some detail.
Denote by $w_1$ the highest weight vector of $\E$, which is assumed to 
be even. Define 
\be
w_i=f_{i-1} w_{i-1}, &\quad 1<i\le n\\
w_0=f_n w_n, &\quad w_{-n} = f_n w_n, \\  
w_{-j}= f_j w_{-j-1}, &\quad n>j\ge 1.
\ee 
Then $\{w_\mu\, |\, \mu=0, \pm 1,\pm 2, ...,\pm n\}$ 
forms a weight basis of $\E$. We will denote by $t$ the 
{\irrep } relative to this basis. The matrix elements of the 
$e_i$, $f_i$ and $k_i$ can be immediately written down. 
We have 
\be 
t(e_i)_{\mu \nu} &=& \delta_{\mu i}\,  \delta_{\nu,  i+1}
                  + \delta_{\mu, -i-1}\,  \delta_{\nu,  -i}, \\
t(f_i)_{\mu \nu} &=& \delta_{\mu,  i+1}\,  \delta_{\nu i}
                  + \delta_{\mu, -i}\,  \delta_{\nu,  -i-1},  \quad i<n, \\
t(e_n)_{\mu \nu} &=& \delta_{\mu n}\,  \delta_{\nu 0}
                  - \delta_{\mu 0}\,  \delta_{\nu,  -n}, \\  
t(f_n)_{\mu \nu} &=& \delta_{\mu 0}\,  \delta_{\nu n}
                  + \delta_{\mu, -n}\,  \delta_{\nu 0},\\
t(k_j)_{\mu \nu} &=& \delta_{\mu \nu}\,  
                q^{(\alpha_j, \, \epsilon_\mu)},  \quad  \quad 1\le j\le n,
\ee  
where $\epsilon_0=0$, and $\epsilon_{-i}=-\epsilon_i$.

Let $\{w^*_\mu\}$  
be the basis of $\E^*$ defined by $w^*_\mu(w_\nu)=\delta_{\mu\nu}$.
$\E^*$ has a natural $\Uq$ - module structure with 
the $\Uq$ action given by 
\bea x w^*_\mu &=&\sum_\nu  (-1)^{[x]\delta_{\mu 0}} 
     t(S(x))_{\mu \nu} w^*_\nu.  \label{duality}\eea   
The lowest weight of $\E$ is $-\epsilon_1$. 
Thus the module $\E$ is self dual.  
This implies that there exists a  
$\Uq$ - module isomorphism $M: \E \rightarrow \E^*$, which 
is unique up to scaler multiples. 
The $w^*_{-1}$, being the highest weight vector of $\E^*$, 
will be identified with $w_1$ so that this arbitrariness in 
$M$ can be removed. Now let  
\be 
w^*_\mu &=&\sum_\nu w_\nu M_{\nu \mu}.  
\ee
Then 
\bea
M_{\mu \nu} &=& m_\mu \delta_{\mu+\nu\, 0},\nonumber\\
m_\mu &=&\left\{ \begin{array}{l l}
                  (-q)^{\mu-1},  &\mu>0,\\
                  (-q)^n,        &\mu=0,\\
                  (-q)^{2n+\mu}, &\mu<0. 
                 \end{array}\right. \label{M} 
\eea

It follows from earlier discussions that  
repeated tensor products of $\E$ are completely reducible. 
Furthermore, every finite dimensional irreducible $\Uq$ 
module is embedded in some $\E^{\otimes k}$ for at least one $k\ge 0$.

For later use, we consider two classes of $\Z_2$ - graded Hopf 
subalgebras of $\Uq$.  Corresponding to  any subset $\bf\Theta$ of $\Nn$,  
we introduce 
\be 
{\cal S}_k &=&\{ k_i^{\pm 1}, i\in \Nn;\ 
               \ e_j, \ f_j, \ j \in {\bf\Theta}\};\\
{\cal S}_{p}&=& {\cal S}_k \cup \{ e_j,
      j \in \Nn\backslash{\bf\Theta}\}.
\ee
The elements of each set generate a $\Z_2$ - graded Hopf subalgebra of
$\Uq$.  The subalgebra generated by the elements of ${\cal S}_k$ will 
be denoted by $\Uk$, and 
called a reductive subalgebra of $\Uq$, while that generated by 
the elements of ${\cal S}_p$ will be denoted by $\Up$ and called 
a parabolic subalgebra. Note that $\Uk$ is a 
$\Z_2$ - graded Hopf subalgebra of $\Up$. 
If we replace $e_i$ by $f_i$ and vice versa in ${\cal S}_{p}$, 
we obtain  another set, which will generate a 
$\Z_2$ - graded Hopf subalgebra of $\Uq$ having similar 
properties as $\Up$.  Results presented in the remainder of the paper can
also be formulated using such algebras.

Observe that there are two types of 
reductive subalgebras, depending on whether $\bf\Theta$ contains $n$. 
The first type arises when $n\not\in \bf\Theta$, and in this case,  
$\Uk$ is the direct product of 
quantized universal enveloping algebras associated with  
a series of ordinary (i.e., non - graded) Lie algebras of 
type $A$ supplemented by the algebra generated by some 
$k_i^{\pm 1}$.  The second type arises when 
$n\in\bf\Theta$.  This time, $\Uk$ is the direct product of 
the first type with a  $\mbox{U}_q(osp(1|2m))$ for some $m<n$.
In both cases, the finite dimensional representations of 
$\Uk$ are completely reducible. 
This fact will be of great importance to the main 
subject of the paper.

Let $V_\mu$ be a finite dimensional irreducible $\Uk$-module.  Then
$V_\mu$ is of highest weight type.  Let $\mu$ be the highest weight
and $\tilde\mu$ the lowest weight of $V_\mu$ respectively.  We can
extend $V_\mu$ in a unique fashion to a $\Up$-module, which is still
denoted by $V_\mu$, such that the elements of ${\cal S}_{p}\backslash
{\cal S}_k$ act by zero.  It is not difficult to see that all finite
dimensional irreducible $\Up$-modules are of this kind.

Consider a finite dimensional irreducible $\Uq$-module
$W(\lambda)$, with highest weight $\lambda$  and lowest weight
$\bar\lambda$.  $W(\lambda)$ can be restricted  in a natural way
to a $\Up$-module, which
is always indecomposable, but not irreducible in general.
It can be readily shown  that
\be
\dim_{\C} \mbox{Hom}_{\Up} ( W(\lambda), \ V_\mu) &=&
    \left\{ \begin{array}{l l}
           1, &  \bar\lambda={\tilde\mu}, \\
           0, &  \bar\lambda\ne {\tilde\mu}.
     \end{array}\right.  \label{Hom}
\ee

\section{\normalsize {\bf THE QUANTUM SUPERGROUP} $\mbox{OSP}_q(1|2n)$}
There exist well established methods  
for quantizing ordinary Lie groups in  
the nonsupersymmetric setting (See \cite{Chari} and references therein.). 
These methods can also be extended 
to construct $\OSPq$, and this will be done here. However, 
we should point out that it is in general much more difficult 
to study quantum supergroups.  See \cite{Zhang} for details 
on $\mbox{GL}_q(m|n)$.  
   
We will show that the quantum supergroup $\OSPq$ admits a 
quantum Haar functional, and also  a Peter - Weyl basis. 
This, however, is an exception rather than the rule. 
It is known that the finite dimensional representations of 
all the quantum superalgebras but $\mbox{U}_q(\osp)$ are 
not completely reducible.  This fact renders it 
impossible to construct Peter - Weyl bases for 
the corresponding quantum supergroups (which are yet to be 
defined except $\mbox{GL}_q(m|n)$.).  
It also forbids the existence of quantum Haar functionals 
on the quantum supergroups in view of the dual Maschke theorem.

Let us recall some general results about $\Z_2$ - graded Hopf algebras. 
Let $A$ be a $\Z_2$ - graded Hopf algebra with co - multiplication 
$\Delta$, co - unit $\epsilon$ and antipode $S$. We define the finite 
dual $A^0$ of $A$ to be a subspace of $A^*$ such that 
for any $f\in A^0$, $Ker f$ contains a two-sided ideal $\cal I$
of $A$ which is of finite co-dimension, i.e.,
$\dim A/{\cal I}<\infty$. Of course in the most general 
situation, there is no guarantee that $A^0$ will not be zero. 
But when $A^0$ is nontrivial, then it is   
also a $\Z_2$ - graded Hopf algebra with a structure 
dualizing that of $A$. More explicitly, the multiplication 
is defined, for $f,\, g\in A^0$, $a,\, b\in A$,  by 
\be 
\langle f g,\ a\rangle &=& \langle f\otimes g, \ \Delta(a)\rangle\\
                &=&\sum_{(a)} (-1)^{[g][a_{(1)}]}
 \langle f,\ a_{(1)} \rangle \langle g, \ a_{(2)} \rangle.
\ee   
It is easy to see that the unit of $A^0$ is $\epsilon$. 
Denote the co - multiplication, the co - unit and the 
antipode of $A^0$ respectively by $\Delta_0$, $\epsilon^0$ and $S_0$. 
Then 
\be
\langle \Delta_0(f),\  a\otimes b\rangle &=&
        \sum_{(f)} (-1)^{[f_{(1)}][f_{(2)}]} \langle f_{(1)}, \ a\rangle 
                   \langle f_{(2)},\ b\rangle \\   
&=&\langle f, \  a b\rangle,\\  
\langle S_0(f),\ a\rangle &=& \langle f,\ S(a)\rangle, \\
\epsilon^0(f) &=& \langle f,\ \1_A \rangle.
\ee

Now we come back to the quantum supergroup $\mbox{OSP}_q(1|2n)$.
As is well known, we can not define the quantum supergroup 
directly. Instead, we need to find the algebra of functions 
on it. Introduce $t_{\mu \nu}\in (\Uq)^*$, 
$\mu, \ \nu=0, \pm 1, \pm 2, ..., \pm n$, defined by 
\be 
t_{\mu \nu}(x)&=& t(x)_{\mu \nu}, \quad \forall x\in\Uq,
\ee
where $t$ is the vector representation of $\Uq$. We call 
the $t_{\mu \nu}$ the matrix elements of $t$. 
Finite dimensionality of $\E$ implies that $t_{\mu \nu}\in(\Uq)^0$,
$\forall \mu,\ \nu$.  

We define the superalgebra $\Tq$ of functions 
on $\OSPq$  to be the 
$\Z_2$-graded  subalgebra 
of $(\Uq)^0$ generated by the matrix elements of the vector
representation of $\Uq$, i.e., $t_{\mu \nu}$, 
$\mu, \ \nu=0, \pm 1, \pm 2, ..., \pm n$.  Then 
\begin{theorem}
\begin{enumerate}
\item $\Tq$ is a $\Z_2$ - graded Hopf algebra.
\item Let $t^{(\lambda)}$ be the \irrep \  of $\Uq$ with 
highest weight $\lambda\in\cP_+$, and let $t^{(\lambda)}_{i j}$,
$i,\, j=1, 2, ..., d_\lambda$ ($d_\lambda= \mbox{dim}t^{(\lambda)})$, 
be the matrix 
elements of $t^{(\lambda)}$. Then 
\bea 
\Tq&=&\bigoplus_{\lambda\in\cP_+} \bigoplus_{i, j=1}^{d_\lambda}
 \C t^{(\lambda)}_{i j}.
\eea 
\end{enumerate}
\end{theorem}  
\begin{proof} {Proof} 
The $\Z_2$ - graded bi - algebra structure of 
$\Tq$ is obvious, and the existence of the antipode 
follows from the self duality of the vector module $\E$ over $\Uq$. 
Part 2) immediately follows from the complete 
reducibility of finite dimensional 
representations of $\Uq$. 
\end{proof} 

Let us now work out the explicit forms of the 
co - multiplication and the antipode.  The co - multiplication is
given by   
\be 
    \Delta_0(t_{\mu \nu}) &=& \sum_{\sigma} (-1)^{(\delta_{\mu 0} 
+\delta_{\sigma 0})(\delta_{\nu 0} + \delta_{\sigma 0})} 
     t_{\mu \sigma}\otimes t_{\sigma \nu}.  
\ee  
The antipode can be constructed from (\ref{duality}) by using the 
$\Uq$ - module isomorphism $M$. We have 
\be 
S_0(t_{\mu \nu})&=& (-1)^{(\delta_{\mu 0} + \delta_{\nu 0})\delta_{\mu 0}} 
                     (M^{-1} t M)_{\nu \mu}\\
                &=& (-1)^{(\delta_{\mu 0} + \delta_{\nu 0})\delta_{\mu 0}} 
                 \ { {m_{-\mu}\ t_{-\nu,\,  -\mu} }\over{m_{-\nu}} },    
\ee
where $m_\mu$ is given by (\ref{M}).

Here we introduce more notations for later use. 
Let $\{w^{(\lambda)}_i | i=1, 2, ..., d_{\lambda}\}$ be 
the homogeneous basis of $W(\lambda)$ with respect to which 
the representation $t^{(\lambda)}$
is defined.  We denote by $\{{\tilde w}^{(\lambda)}_i 
| i=1, 2, ..., d_{\lambda}\}$
the basis of $W(\lambda)^*=W(\lambda^\dagger)$ such that
${\tilde w}^{(\lambda)}_i(w_j)=\delta_{i j}.$
The $\Uq$ - module structure of $W(\lambda)^*$ 
enables us to define 
${\tilde t}^{(\lambda)}_{i j}\in\Tq$ by
\ban
x {\tilde w}^{(\lambda)}_i &=&\sum_{j} {\tilde t}^{(\lambda)}_{j i}(x) 
{\tilde w}^{(\lambda)}_j,
\quad \forall x\in\Uq.
\nan
Then
\ban
{\tilde t}^{(\lambda)}_{j i} 
&=& (-1)^{[i]([i]+[j])} S_0(t^{(\lambda)}_{i j}),
\nan
where $[i]=0$ or $1$ depending on whether $w_i$ is even or odd.
Cleary the ${\tilde t}^{(\lambda)}_{j i}$ are linear combinations
of $t^{(\lambda^\dagger)}_{i j}$. Furthermore, 
the ${\tilde t}^{(\lambda)}_{j i}$, $\forall\lambda\in{\cal P}_+$,  
also form a basis of $\Tq$.

From here on, we will omit the subscript $_0$ from $\Delta_0$ 
and $S_0$.

Let us now turn to the discussion of 
a Haar functional on the quantum supergroup $\Tq$.
But before embarking on this task, we first consider the notion 
of an integral on an arbitrary  $\Z_2$ - graded Hopf algebra $A$. 
Let $A^*$ be its dual, which has a natural
$\Z_2$ - graded algebraic structure induced by the co-algebraic
structure of $A$. 
An even homogeneous element $\int^l\in A^*$ is called a
left integral on $A$ if
\ban
f\cdot \int^l &=& \langle f, \1_A\rangle\ \int^l, \ \ \ \forall f\in A^*.
\nan
Similarly,  an even homogeneous element 
$\int^r\in A^*$ is called a right integral on $A$ if
\ban
\int^r\cdot f &=& \langle f, \1_A\rangle \ \int^r, \ \ \ \forall f\in A^*.
\nan
A straightforward calculation shows that the defining properties of
the integrals are equivalent to the following requirements
\ba
(\id\otimes \int^l)\Delta(x)=\int^l x, \quad
(\int^r\otimes \id)\Delta(x)=\int^r x,\quad  \forall x\in A.
\na
where ${\rm id}$ is the identity map on $A$.

A Haar functional $\int\in A^*$ on $A$ is an integral on $A$ which is
both left and right, and sends $\1_A$ to $1$, i.e.,
\ba
&(i).&  (\int\otimes \id) \Delta(x) = (\id\otimes \int)\Delta(x)
       = \int x, \quad \forall x\in A,\nonumber\\
&(ii).&  \int \1_A  = 1. \label{integral}
\na

In the case of $\Tq$,  it is an entirely straightforward 
matter to show that  
\begin{theorem}
The element $\int\in (\cT)^*$ defined by
\ba
\int \1_{\cT}=1; \quad  \quad \int t^{(\lambda)}_{i j}=0,
\quad  0\ne \lambda\in \cP_+, \nonumber
\label{Haar}
\na
gives rise to a Haar functional on $\cT$.
\end{theorem}

Denote by $2\rho$ the sum of the positive roots of $\frak g$.
Let $K_{2\rho}$ be the product of powers of $k_i^{\pm 1}$'s such that
\ban 
K_{2\rho} e_i K_{2\rho}^{-1} &=& q^{(2\rho, \ \alpha_i)} e_i,
\quad \forall i.\nan 
Then it can be easily shown  that 
\ban 
S^2(x) &=& K_{2\rho} x K_{2\rho}^{-1}, \quad \forall x\in\Uq.\nan 
We define the quantum super-dimension of the  
irreducible $\Uq$-module  $W(\lambda)$ by 
\ban 
\mbox{SD}_q(\lambda)&:=& Str\{t^{(\lambda)}(K_{2\rho})\}.   
\nan 
For quantum superalgebras other than the $\osp$ series, 
there exists a class of finite dimensional irreducible representations, 
the typicals, of which the super-dimensions vanish identically.  
Again, $\mbox{U}_q(\osp)$ is an exception, and we have the following 
important property: for any irreducible $\mbox{U}_q(\osp)$ - 
module $W(\lambda)$ with highest weight $\lambda\in\cP_+$,  
\ban
\mbox{SD}_q(\lambda)&\ne& 0.  
\nan 
Now the Haar functional $\int$ satisfies the following properties.
\begin{lemma}\label{tresult}
\ba 
\int t^{(\lambda)}_{i j} {\tilde t}^{(\mu)}_{r s} 
(-1)^{[j] [r] + [i]+[j]}
&=& \delta_{i r} \delta_{\lambda  \mu} 
{{t^{(\lambda)}_{s j}(K_{2\rho})}\over{\SDq(\lambda)} }, \nonumber \\
\int {\tilde t}^{(\lambda)}_{i j} t^{(\mu)}_{r s} 
(-1)^{[j][r]}&=&
\delta_{j s} \delta_{\lambda  \mu} 
{{{\tilde t}^{(\lambda)}_{ i r}(K_{2\rho})}\over{\SDq(\lambda)} }.
\na 
\end{lemma}
\begin{proof} {Proof}
Consider the first equation.  
The $\lambda \ne \mu$ case is easy to prove:
the integral vanishes because
the tensor product $W(\lambda)\otimes W(\mu^\dagger)$
does not contain the trivial $\Uq$-module.
When   $\lambda = \mu$,  we introduce the notations
\ban
\phi_{i r; s j} = \int t^{(\lambda)}_{i j}
                    {\tilde t}^{(\lambda)}_{r s} 
                  (-1)^{[j] [r] + [i]+[j]};
&\quad \Phi[s, j] = \left( \phi_{i r; s j}\right)_{i, r=1}^{d_\lambda};
&\quad \Psi[i, r]= \left( \phi_{i r; s j}\right)_{s, j=1}^{d_\lambda}.
\nan
It is clearly true that
$Str\left(\Psi[i, r]\right)$ $=$ $\delta_{i r}$.
 
Note that corresponding to each $x\in\Uq$, there exists an
${\tilde x}\in (\cT)^*$ defined by
${\tilde x}(a) = \langle a, \ x\rangle$, $\forall a\in\cT$.
The left integral property of $\int$ leads to
\ban
\epsilon(x)\phi_{i r; s j} &=& ({\tilde x}.\int)t^{(\lambda)}_{i j}
{\tilde t}^{(\lambda)}_{r s} (-1)^{[j] [r] + [i]+[j]}\\
&=& \sum_{(x)} \sum_{i', r'} t^{(\lambda)}_{i i'}(x_{(1)})
t^{(\lambda)}_{ r' r}(S(x_{(2)})) \phi_{i' r'; s j} (-1)^{[x]([i]+[j]) 
+ [x_{(2)}]([j]+[s])},\\
\mbox{i.e.}\quad  \epsilon(x) \Phi[s, j]
&=& \sum_{(x)} t^{(\lambda)}(x_{(1)}) \Phi[s, j] t^{(\lambda)}(S(x_{(2)}))
(-1)^{[x_{(2)}]([j]+[s])},
\quad \forall x\in\Uq. 
\nan
Schur's lemma forces $\Phi[s, j]$ to be proportional to the identity 
matrix, and we have
\ban
\Psi[i, r]&=&\delta_{i r} \psi,
\nan
for some $d_\lambda \times d_\lambda$ matrix $\psi$.
The right integral property of $\int$ leads to
\ban
\epsilon(y) \psi &=& \sum_{(y)} t^{(\lambda)}(K_{2\rho})
t^{(\lambda)}(y_{(1)}) t^{(\lambda)}(K_{2\rho}^{-1}) \psi
t^{(\lambda)}(S(y_{(2)})).
\nan
Again by using Schur's lemma we conclude that $\psi$ is proportional
to $t^{(\lambda)}(K_{2\rho})$. Since its supertrace is $1$, we have
\ban \psi &=& { {t^{(\lambda)}(K_{2\rho})}\over{\SDq(\lambda)} }. \nan
This completes the proof of the first equation of the lemma.
The second equation can be shown in exactly the same way.
\end{proof}

It is worth observing that this Lemma and part 2) of Theorem 1 
provide a quantum analog of the Peter - Weyl 
theorem for $\OSPq$.

\section{\normalsize{\bf QUANTUM HOMOGENEOUS SUPERVECTOR BUNDLES}}
In this section we will investigate the quantum homogeneous 
superspaces and quantum homogeneous supervector bundles 
arising from the quantum supergroup $\OSPq$ by adapting the 
methods and techniques of \cite{Zhang, Gover} to the 
present context. 
Let us start by introducing two types of actions of $\Uq$ on $\cT$. 
The first action will be denoted by $\circ$, which 
corresponds to the right translation in the classical theory 
of Lie groups.  It is defined by 
\ba
x\circ f&=&\sum_{(f)} (-1)^{[f_{(1)}][f_{(2)}]} \ 
         f_{(1)} \ \langle f_{(2)}, \ x\rangle, 
\quad x\in\Uq, \ f\in\cT,  \label{circ}
\na
Straightforward calculations show that
\ban x\circ (y \circ f ) &=& (x y)\circ f;\\
(x\circ f) ( y ) &=& f(y x), \\
(id_{\cT}\otimes x\circ) \Delta(f) &=& \Delta( x\circ f ).  \nan
The other action, which corresponds to the
left translation in the classical  Lie group theory,
will be denoted by $\cdot$. It is defined by  
\ba
x\cdot f &=&\sum_{(f)} \langle f_{(1)}, \ S^{-1}(x) \rangle f_{(2)}. 
\na
It can be easily shown that  
\ban
(x\cdot f)(y)&=& (-1)^{[x][y]} f( S^{-1}(x) y ),\\
x\cdot(y\cdot f)&=& (x y)\cdot f, \quad  x,\ y\in \Uq, \ f\in \cT.
\nan
Furthermore, the two actions graded - commute in the following sense  
\ban 
x\circ(y\cdot f)&=&(-1)^{[x][y]} y\cdot(x\circ f), 
\quad  x, y\in\Uq, \ f\in\cT.
\nan 

Let $V$ be a finite dimensional module over $\Uk$.
We extend the actions $\circ$ and $\cdot$ trivially to 
$V\otimes \Tq$: for any              
$\zeta= \sum v_i \otimes f_i  \in V\otimes\Tq$,
\ban
x\cdot\zeta &=& \sum (-1)^{[x][v_i]} v_i \otimes x\cdot f_i, \\
x\circ\zeta&=& \sum (-1)^{[x][v_i]} v_i \otimes x\circ f_i,
\ \ \ \  x\in\Uq.
\nan

We now introduce two important definitions:
\ba
\Eq&:=&\left\{ f\in\Tq \ | \ x\circ f =\epsilon(x) f, 
\quad \forall x\in\Uk \right\};   \label{space}\\  
\Hq(V)&:=& \left\{ \zeta\in V\otimes \Tq \
| \ x\circ\zeta = (S(x)\otimes id_{\ca}) \zeta, \ \forall x\in\Uk\right\}.
\label{section}
\na

The remainder of this section is devoted to studying the properties 
of these objects. Let us first prove the following   
\begin{proposition}\label{plenty}
1). $\Eq$ is an infinite dimensional subalgebra of $\Tq$. \\
2). $\Hq(V)$ is an infinite dimensional supervector space 
if the weight of any vector of $V$ is $\Uq$-integral, and is 
zero otherwise. 
\end{proposition}
\begin{proof} {Proof}
We first show that 
$\Eq$ is a subalgebra of $\Tq$. 
Since $\Uk$ is  a Hopf subalgebra of $\Uq$, for any $x\in\Uk$, 
$\Delta(x)=\sum_{(x)}x_{(1)}\otimes x_{(2)}$ $\in\Uq\otimes\Uq$. 
Hence 
\ban 
x\circ ( a b )&=&
\sum_{(x)} (-1)^{[x_{(2)}][a]}  \{x_{(1)}\circ a\} \{ x_{(2)}\circ b\} 
   =\epsilon(x) ab,  
\nan 
that is, $a b\in\Eq$.

Since the finite dimensional representations of $\Uk$ are completely 
reducible, the study of properties of $\Hq(V)$  reduces to the case 
when $V$ is irreducible. Let $V_\mu$ be a
finite dimensional irreducible $\Uk$ -module with highest weight $\mu$
and lowest weight $\tilde\mu$. 
Any element $\zeta$ $\in$  $\Hq(V_\mu)$
can be expressed in  the form
\ban
\zeta &=& \sum_{\lambda\in\cP_+}
\sum_{i, j} v_{i j}^{(\lambda)}\otimes {\tilde t}_{i j}^{(\lambda)}, 
\nan
for some $v_{ i j}^{(\lambda)}\in V_\mu$.
Fix an arbitrary $\lambda\in\cP_+$. For any nonvanishing
$w\in W(\lambda)$, the following linear map is clearly surjective: 
\ban
\mbox{Hom}_{\C}(W(\lambda), \ V_\mu)\otimes w &\rightarrow& V_\mu,\\
\phi\otimes w&\mapsto& \phi(w). 
\nan
Thus there exist $\phi^{(\lambda)}_i$
$\in$ $\mbox{Hom}_{\C}(W(\lambda), \ V_\mu)$ such that
$v_{ i j}^{(\lambda)} = \phi^{(\lambda)}_i ( w_j^{(\lambda)} )$,
where $\{w_i^{(\lambda)}\}$ is the basis of $W(\lambda)$ discussed 
before. Therefore,  we can rewrite $\zeta$ as
\ban
\zeta &=& \sum_{\lambda\in\cP_+}
\sum_{i, j} \phi_i^{(\lambda)}(w_j^{(\lambda)})
\otimes {\tilde t}_{i j}^{(\lambda)}.
\nan
The defining property of $\Hq(V_\mu)$ states that
\ban
\ell\circ\zeta&=& ( id_{\cT}\otimes S(\ell)) \zeta, \quad \forall \ell\in\Uk .
\nan
Thus  we have
\ban
\sum_{\lambda\in\cP_+} \sum_{i, j, k} 
t_{j k }^{(\lambda)}(S(\ell)) \phi_i^{(\lambda)}(w_j^{(\lambda)}) 
\otimes (-1)^{[\ell] [\phi_i^{(\lambda)}]} {\tilde t}_{i k}^{(\lambda)}
&=& \sum_{\lambda\in\cP_+} \sum_{i, j} S(\ell)
\phi_i^{(\lambda)}(w_j^{(\lambda)}) \otimes {\tilde t}_{i j}^{(\lambda)}.
\nan 
Recalling that the ${\tilde t}_{k i}^{(\lambda)}$ are linearly
independent. So the above  is equivalent to
$$
\ell \phi_i^{(\lambda)}(w_j^\lambda)
= (-1)^{[\ell] [\phi_i^{(\lambda)}]}
\phi_i^{(\lambda)}(\ell w_j^\lambda), ~~~ \forall \ell\in\Uk .
$$
This equation is precisely the statement that the $\phi_i^{(\lambda)}$
be $\Uk$-module homomorphisms of degrees $[\phi_i^{(\lambda)}]$,
\ban
\phi_i^{(\lambda)}&\in&\mbox{Hom}_{\Uk}\left(W(\lambda), V_\mu\right)
\subset \mbox{Hom}_{\C}\left(W(\lambda), V_\mu\right),
\quad \forall i.
\nan
Thus finding sections in $\Hq(V_\mu)$ is equivalent to finding, 
for all $\lambda\in\cP_+$, the homomorphisms $\phi^{(\lambda)}$ $\in$ 
$\mbox{Hom}_{\Uk}\left(W(\lambda), V_\mu\right)$. 
Note that each such homomorphism $\phi^{(\lambda)}$
determines $d_\lambda$ linearly independent sections
$$
\zeta^{(\lambda)}_i =  \sum_{j}  
\phi^{(\lambda)}(w_j^{(\lambda)})\otimes {\tilde t}_{i j}^{(\lambda)}. 
$$

However, when $\mu$ is not integral with respect to $\Uq$, 
$\mbox{Hom}_{\Uk}\left(W(\lambda), V_\mu\right)=0$, and hence 
$\Hq(V_\mu)$ vanishes in this case. 

Now consider the case with $\mu=0$, we have $\Hq(V_{\mu=0})=\Eq$ 
as supervector spaces. 
There is a homomorphism from the trivial representation of $\Uq$,
$W(0)={\Bbb C}$, onto $V_0={\Bbb C}$. This gives the constant sections of
$\Eq$.  Let $\gamma$ be the highest root of $\frak g$. 
Recall that in the classical situation, 
${\frak k}$ is reductive with  
$N=r-|{\bf\Theta}|$ independent central elements. 
This,  transcribed to the quantum case, 
implies the existence of $N$ linearly independent $\Uk$-homomorphisms
$W(\gamma)\to {\Bbb C}$. 
As mentioned above each of these corresponds
to $d = \mbox{dim} ({\frak g})$  linearly independent sections.  
So the representation
$W(\gamma)$ determines $ N d$ 
linearly independent sections.
Further linearly independent sections can  be obtained using the 
following lemma
\begin{lemma} \label{prodhom}
Suppose there are non-trivial $\Uk$-homomorphisms
$W(\lambda_1)\to V_{\mu_1}$
and $W(\lambda_2)\to V_{\mu_2}$. Then there is an induced non-trivial
$\Uk$-homomorphism
$$
W(\lambda_1+\lambda_2)\to V_{\mu_1+\mu_2}.
$$
\end{lemma}
For example, for any positive integer $m$, there exist $(m|N)$ 
(partition of $m$ into $\leq N$ parts) linearly independent 
homomorphisms $W(m\gamma)\to {\Bbb C}$.  Thus we have proved that  
the algebra $\Eq$ is infinite dimensional.

Now let us consider the case with $0\neq\mu\in \cP$. 
It is an elementary exercise to verify that $V_\mu$ is 
$\Uk$-isomorphic to a $\Uk$-irreducible part of $W(\lambda')$, 
where $\lambda'$ is the dominant weight in the Weyl group orbit 
of $\mu$. Thus there is a non-trivial $\Uk$-homomorphism 
$$ W(\lambda') \to V_\mu , $$
and this  determines at least  $d_{\lambda'}$ linearly independent 
sections in $\Hq(V_\mu)$. Further linearly independent sections 
can be constructed explicitly using lemma
\ref{prodhom} which promises a family of homomorphisms
$$
W(\lambda'+m\gamma)\to V_\mu ~~~~~m\in {\Bbb N}_+.
$$
This establishes that $\Hq(V_\mu)$ is infinite dimensional. 
\end{proof}

$\Eq$ may be regarded as the quantum analog of the algebra of 
functions over the superspace $\mbox{OSP}(1|2n)/K$, where $K$ 
is the subgroup of $\mbox{OSP}(1|2n)$ with Lie superalgebra ${\frak k}$. 
Such homogeneous superspaces were studied in the work of Manin
\cite{Manin}, Penkov \cite{Penkov} and others.  
Here we wish to make some investigations into their quantum analogs.

As is well known,  one can not define a noncommutative (in the 
$\Z_2$ - graded sense) space directly in geometrical 
terms. Instead,  such a space  has to be defined by specifying 
its algebra of functions.  We will take $\Eq$ as 
the algebra of functions over the quantum homogeneous superspace    
which corresponds to $\mbox{OSP}(1|2n)/K$ in the classical situation. 
Let us now study properties of $\Hq(V)$. First observe that  
\begin{theorem}
$\Hq(V)$ furnishes  a two-sided  $\Eq$ module 
     under the multiplication of $\Tq$.
\end{theorem}
\begin{proof} {Proof}  
The left and right actions of $\Eq$ 
on $\Hq(V)$ are respectively defined by 
\ban 
a\zeta&=&\sum_r (-1)^{[a][v_i]} v_i\otimes a f_i,\\
\zeta a&=&\sum_r v_i\otimes f_i a,  
\nan  
where $a\in\Eq$ and $\zeta=$ $\sum_i v_i\otimes f_i$ $\in\Hq(V)$. 
Now for $p\in\Uk$,  
\ban 
p\circ(a\zeta) &=& \sum_{(p)} (-1)^{[p_{(2)}][a]} 
                 \{ p_{(1)}\circ a \} \{p_{(2)}\circ \zeta\}\\
               &=& (-1)^{[p][a]} a\{p\circ \zeta\} 
                = (S(p)\otimes id_{\ca}) a\zeta;\\
p\circ(\zeta a)&=& \sum_{(p)} (-1)^{[p_{(2)}][\zeta]} 
                \{p_{(1)}\circ \zeta\} \{ p_{(2)}\circ a \}\\
               &=&\{p\circ \zeta\}a= (S(p)\otimes id_{\ca}) \zeta a.
\nan 
This completes the proof.
\end{proof} 

When $V$ is actually a $\Uq$ - module, 
the $\Eq$ - module $\Hq(V)$ has a particularly simple structure.
\begin{proposition}\label{projective}
Let $W$ be a finite dimensional left $\Uq$-module, which we
regard as a left $\Uk$-module by restriction.
Then $\Hq(W)$ is isomorphic to
$W\otimes \Eq$ either as a left or right $\Eq$-module.
\end{proposition}
\begin{proof} {Proof}
We first construct the right $\Eq$-module 
isomorphism. 
Being a  left $\Uq$-module, $W$ carries a natural 
right $\ca$  co-module structure with the co-module action 
$\delta: W \rightarrow W\otimes\cT$ 
defined by 
\ba 
\delta(w)(x) &=& x w, \quad  x\in\Uq, \ \ w\in W. 
\label{comodule} 
\na
(Here the notation requires some clarification. If we express 
$\delta(w)=\sum_{(w)} w_{(1)}\otimes w_{(2)}$, 
then $\delta(w)(x) = \sum_{(w)} (-1)^{[x][w_{(1)}]} 
w_{(1)} \langle w_{(2)}, \, x\rangle$.) 
Define $\eta: W\otimes\Tq \longrightarrow W\otimes\Tq$ 
by the composition of maps
\ban 
W\otimes\Tq  \stackrel{\delta\otimes id}{\longrightarrow} 
W\otimes \Tq\otimes\cT  
\longrightarrow W\otimes\Tq, 
\nan   
where the last map is the multiplication of $\ca$.
Then $\eta$ defines a right $\Eq$-module isomorphism, 
with the inverse map given by the composition 
\ban
W\otimes\Tq \stackrel{\delta\otimes id}{\longrightarrow}
W\otimes\Tq\otimes\Tq  
\stackrel{(id\otimes S\otimes id)}{\longrightarrow}
W\otimes \cT\otimes \Tq\longrightarrow W\otimes \Tq, 
\nan
where the last map is again  the multiplication of $\ca$.
It is not difficult to show that 
\ban 
x\circ\eta(\zeta)&=&\eta\left(\sum_{(x)}(x_{(1)}\otimes  id_{\ca}) 
                     x_{(2)}\circ \zeta\right),\\
x\circ\eta^{-1}(\zeta)&=&\eta^{-1}\left(\sum_{(x)} 
(S(x_{(1)})\otimes  id_{\ca}) x_{(2)}\circ\zeta\right), 
\quad \forall \zeta\in\Tq\otimes W, \ \  x\in\Uq.
\nan 
Consider $\zeta\in\Hq(W)$. We have 
\ban 
p\circ\eta(\zeta)&=&  
  \eta\left(\sum_{(p)} (p_{(1)}\otimes  
   id_{\ca})p_{(2)}\circ\zeta\right)\\
&=&\eta\left(\sum_{(p)}  (p_{(1)} S(p_{(2)})\otimes  
   id_{\ca})\zeta\right)\\ 
&=&\epsilon(p)\eta(\zeta), \quad  \forall \ p\in\Uk.
\nan 
Hence $\eta(\Hq(W))\subset W\otimes\Eq$. 
Conversely, given any $\xi\in  W\otimes \Eq$, we have 
\ban 
p\circ\eta^{-1}(\xi)&=& 
\eta^{-1}\left(\sum_{(p)} (S(p_{(1)})\otimes  id_{\ca}) 
p_{(2)}\circ \xi\right)\\
&=&\eta^{-1}\left(\sum_{(p)}(S(p_{(1)})\epsilon(p_{(2)}) 
\otimes  id_{\ca})\xi\right)\\   
&=& (S(p)\otimes  id_{\ca}) \eta^{-1}(\xi),  \quad  \forall \ p\in\Uk.   
\nan 
Thus $\eta^{-1}(W\otimes \Eq)\subset\Hq(W)$. 
Therefore the restriction of $\eta$ to $\Hq(W)$  provides 
the desired right $\Eq$-module isomorphism. 

The left module isomorphism is given by the restriction to $\Hq(W)$ 
of the linear map  
$\kappa: W\otimes \Tq$ $\rightarrow$ 
$W\otimes \Tq$, which is defined by the following composition of maps 
\ban
W\otimes \Tq \stackrel{\delta\otimes id}{\longrightarrow}
W\otimes\Tq\otimes\Tq 
\stackrel{id\otimes P(S^2\otimes id)}{\longrightarrow}
W\otimes\Tq\otimes \Tq\longrightarrow W\otimes \Tq,
\nan
where 
\ba 
P: \Tq\otimes\Tq&\rightarrow &\Tq\otimes \Tq, \nonumber \\
   a\otimes b &\mapsto& (-1)^{[a][b]}\,  b \otimes a. \label{P}  
\na 
The inverse map $\kappa^{-1}$ is given by 
\ban
W\otimes \Tq \stackrel{\delta\otimes id}{\longrightarrow}
W\otimes\Tq\otimes\Tq 
\stackrel{id\otimes P(S\otimes id)}{\longrightarrow}
W\otimes\Tq\otimes \Tq\longrightarrow W\otimes \Tq.  
\nan
\end{proof} 

With the help of this Proposition, we can now prove the following 
important result. 
\begin{theorem}\label{proj}
$\Hq(V)$ is projective and of finite type
both as a left and right module over the superalgebra $\Eq$ of
functions on the quantum homogeneous superspace.
\end{theorem}
\begin{proof} {Proof} 
Since $\Uk$ is a reductive subalgebra of $\Uq$,
all finite dimensional representations of $\Uk$ 
are completely reducible. 
Let $V_s$, $s=1, 2, ..., K\le\infty$, be the irreducible 
direct summands of $V$ such that their weights are all integral 
with respect to $\Uq$. Then $\Hq(V)=\bigoplus_s \Hq(V_s)$. 
Consider any $V_s$, and denote its highest weight by $\mu_s$.  
There exists such a ${\hat\mu}_s$ in the Weyl group orbit of 
${\frak g}$ that is integral dominant with respect to $\frak g$. 
Let $W({\hat\mu}_s)$ be the irreducible $\Uq$ - module 
with highest weight ${\hat\mu}_s$, which can be regarded 
as a $\Uk$ - module in the natural way. There always exists a 
$\Uk$ - module $V_s^\bot$ such that $W({\hat\mu}_s) 
= V_s\oplus V_s^\bot$. 
Write $V^\bot=\oplus_s V_s^\bot$, and $W=\oplus_s W({\hat\mu}_s)$. 
We have  
\ban
\Hq(V) \oplus \Hq(V^\bot) &=&\Hq(W)\\
&\cong& W\otimes\Eq, 
\nan
where the last step follows from Proposition \ref{projective}.
\end{proof}

Recall that in classical differential geometry, 
the space $\cal H$ of sections of
a vector bundle over a compact manifold $M$ furnishes a module over
the algebra ${\cal A}(M)$ of functions on $M$.  
It then follows from Swan's
theorem that this module must be projective and is of finite type.
Conversely, any projective module of finite type over ${\cal A}(M)$ 
is isomorphic to the space of sections of some vector bundle over $M$.  
This result is taken
as the starting point for studying vector bundles in noncommutative
geometry: one defines a vector bundle over a noncommutative space in
terms of the space of sections which is required to be a finite type
project module over the noncommutative algebra 
of functions on the virtual noncommutative space.
Therefore, $\Hq(V)$ will be called the space of
sections of a quantum supervector bundle over the 
quantum homogeneous superspace associated with $\Eq$. 

Homogeneous supervector bundles at the classical level
were studied in \cite{Manin, Penkov}. We will not enter 
the discussion of the subject, but refer to 
the original papers for detail. It is worth mentioning 
that the supergeometry of these objects are very interesting, 
and many aspects of it remain to be understood. 
The Hopf algebra approach to quantum homogeneous supervector bundles
adopted here works equally well at the classical level, 
and should provide a new approach to the supergeometry of 
homogeneous supervector bundles.

Follwoing the classical terminology, we will call a 
quantum supervector bundle trivial if the sections form 
a free module over the superalgebra of functions on the 
quantum superspace. It immediately follows from 
Proposition \ref{projective} that 
\begin{proposition}
If the $\Uk$ - module $V$ is in fact a finite dimensional left $\Uq$
- module, then the quantum homogeneous supervector bundle with the 
space of sections $\Hq(V)$ is trivial.
\end{proposition}

\section{\normalsize{\bf INDUCED REPRESENTATIONS}}
In this section we will investigate induced representations of the 
quantum supergroup $\OSPq$ by using results of the last section.
The following Proposition explains how quantum homogeneous 
supervector bundles enter representation theory. 
\begin{proposition}
$\Gamq(V)$ furnishes a left $\Uq$ - module under the $\cdot$ action,
and also a right $\Tq$ co - module under the action 
$\omega=\mbox{id}_V\otimes (\mbox{id}_{\Tq}\otimes S^{-1})\Delta$.  
\end{proposition}
\begin{proof} {Proof}  
For $p\in\Uk$, $x\in\Uq$, and $\zeta\in\Hq(V)$, we have 
\ban
p\circ  ( x\cdot\zeta )&=& (-1)^{[p][x]} x\cdot(p\circ\zeta)\\
&=&(S(p)\otimes id_{\ca}) ( x\cdot\zeta ).   
\nan
Thus $\Gamq(V)$ indeed furnishes a left $\Uq$ - module under 
the $\cdot$ action.  The $\Tq$ - co-action $\omega$ is 
just the dual of this left $\Uq$ -  action.  
\end{proof}

We call $\Hq(V)$ an induced $\Uq$ module, and also 
an induced $\Tq$ co-module.
For such induced modules, we have the following
quantum analog of Frobenius reciprocity.
\begin{theorem} Let $W$ be a $\Uq$ module, 
the restriction of which furnishes a $\Uk$ module in a natural way.
Then there exists a canonical isomorphism
\ba
\mbox{Hom}_{\Uq} ( W, \ \Hq(V)  )
&\cong& \mbox{Hom}_{\Uk} ( W, \  V),  
\na 
where $\Uq$ acts on the left module $\Hq(V)$ via the $\cdot$ action. 
\end{theorem}
\begin{proof} {Proof}
  We prove the Proposition by explicitly constructing
the isomorphism, which we  claim to be  the linear map
\ban
F: \mbox{Hom}_{\Uq} ( W, \ \Hq( V ) ) &\rightarrow&
   \mbox{Hom}_{\Uk} ( W, \ V), \\
    \psi &\mapsto& \psi(1_{\Uq}),    
\nan
with the inverse map  
\ban
\bar{F}: \mbox{Hom}_{\Uk} ( W, \ V) &\rightarrow&
        \mbox{Hom}_{\Uq} ( W, \ \Hq( V ) ), \\
       \phi &\mapsto& \bar{\phi}=(\phi \otimes S) \delta,   
\nan 
where $\delta: W\rightarrow W\otimes\Tq$ is 
the right $\ca$ co-module action defined by (\ref{comodule}).

To verify our claim, we first need to demonstrate that 
the  image of $F$  is contained in
$\mbox{Hom}_{\Uk} ( W, \ V)$.  Consider $\psi$ $\in$ 
$\mbox{Hom}_{\Uq} ( W, \ \Hq( V ) )$.  For any $p$ $\in$ $\Uk$ and 
$w$ $\in$ $W$,  we have 
\ban
p ( F \psi (w) ) &=& (S^{-1}(p)\circ\psi (w)) (\1_{\Uq}),  
\nan
where we have used the defining property of $\Hq( V )$.  
Note that 
\ban 
(S^{-1}(p)\circ\psi (w)) (\1_{\Uq})&=& (p\cdot\psi (w)) (\1_{\Uq}). 
\nan  
The $\Uq$-module structure of $\Hq( V )$ and the given condition 
that $\psi$ is a $\Uq$-module homomorphism immediately leads to 
\ban
p ( F \psi (w) ) &=& (-1)^{[\psi][p]}\,  \psi(p w) (\1_{\Uq}) \\
                 &=&(-1)^{[\psi][p]}\,  F \psi (p w), 
\quad   p\in \Uk, \  w\in W. 
\nan  

Now consider ${\bar F}$.
We first show that the image $Im({\bar F})$
of ${\bar F}$ is contained in
$\mbox{Hom}_{\Uq} ( W, \ \Hq( V ) )$.
Note that $Im({\bar F})$ $\subset$  $\mbox{Hom}_{\C}
(W, \  V\otimes \Tq)$. 
Some relatively simple manipulations lead to
\ban
(x\cdot{\bar\phi} (w) )  &=& {\bar\phi} ( x w ),\\
( p\circ {\bar\phi} (w) ) &=& (S(p)\otimes id_{\ca}) 
                           {\bar\phi}(w), \quad 
x\in\Uq, \ \ p\in\Uk, \ \ w\in W.
\nan
Therefore,  $Im({\bar F})$ $\subset$
$\mbox{Hom}_{\Uq} ( W, \ \Hq( V ) )$.
Now we show that $F$ and ${\bar F}$ are inverse to each other.
For $\psi\in$ $\mbox{Hom}_{\Uq} ( W, \ \Hq( V ) )$, and
$\phi$ $\in$ $\mbox{Hom}_{\Uk} ( W, \ V)$, we have
\ban
(F \bar{F}  \phi)(w) &=& (\bar{F} \phi)(w) (1_{\Uq})\\
                             &=& \phi(w), \\
(\bar{F} F \psi)(w)(x)&=& (-1)^{[x]([w]+1)}\  ( F \psi) (S(x) w)\\
        &=&(-1)^{[x]([w]+1)}\   \psi (S(x) w) (1_{\Uq})\\
        &=&(-1)^{[x]([w]+[\psi]+ 1)}\    (S(x)\cdot \psi (w)) (1_{\Uq})\\
        &=& \psi (w) (x), \ \ \ \ x\in\Uq, \ \ w\in W.
\nan
This completes the proof of the Proposition.
\end{proof} 
 
Let $V_\mu$ be a finite dimensional irreducible $\Up$-module with
highest weight $\mu$ and lowest weight $\tilde\mu$. 
Since $V_\mu $ is a $\Up$-module the following is a well defined 
subspace of $\Hq(V_\mu)$,
\ban \co_q(V_\mu) & :=& \left\{ \zeta\in \Hq(V_\mu) \ | 
\ p\circ \zeta =
(S(p) \otimes id_{\ca})\zeta , \ \forall p\in\Up\right\}.
\nan
We may regard $\co_q(V_\mu)$ 
as the quantum analog of the space of `holomorphic sections'. 
Recall that the notation $W(\lambda)$ denotes
the irreducible $\Uq$ module with highest weight $\lambda$.  We have
the following result.
\begin{theorem} \label{BW}
There exists the following $\Uq$ module isomorphism
\ba
\co_q(V_\mu) &\cong&\left\{\begin{array}{ l l }
       W((-{\tilde \mu})^\dagger),& -{\tilde \mu}\in\cP_+,\\
       0,& \mbox{otherwise}.
       \end{array} \right.
\na
\end{theorem}
\begin{proof}{Proof} 
Each  $\zeta$ $\in$ $\co_q(V_\mu)$ 
can be expressed in the form 
\ban 
\zeta &=& \sum_{\lambda\in\cP_+}
\sum_{i, j} v_{i j}^{(\lambda)}\otimes  {\tilde t}_{i j}^{(\lambda)}, 
\nan
for some $v_{ i j}^{(\lambda)}\in V_\mu$ ($i,j=1,\cdots d_\lambda$).
Arguing as in the proof of proposition \ref{plenty}
one concludes, for each $\lambda\in\cP_+$, that
there exist $\phi^{(\lambda)}_i$ $\in$ $\mbox{Hom}_{\C}(W(\lambda), \ 
V_\mu)$ such that $v_{ i j}^{(\lambda)} = \phi^{(\lambda)}_i (
w_j^{(\lambda)} )$, where $\{w_i^{(\lambda)}\}$ is the basis of
$W(\lambda)$, relative to which the \irrep \  $t^{(\lambda)}$ of $\Uq$ is
defined.  Thus we can rewrite $\zeta$ as 
\ban
\zeta &=& \sum_{\lambda\in\cP_+} \sum_{i, j} 
  \phi_i^{(\lambda)}(w_j^{(\lambda)})
\otimes {\tilde t}_{i j}^{(\lambda)}.  
\nan
Similar reasoning as in the proof of proposition \ref{plenty}
shows that 
the $\phi_i^{(\lambda)}$ must be $\Up$-module homomorphisms of
degree $[\phi_i^{(\lambda)}]$. 
It immediately follows from (\ref{Hom}) that  
\ban  \phi^{(\lambda)}_i &=& c_i \, \phi^{(\lambda)},\quad 
c_i\in \C, \nan  
and $\phi^{(\lambda)}$ may be nonzero only when
\ban
{\bar \lambda}&=&\tilde\mu. 
\nan
Hence, if $-\tilde\mu\not\in\cP_+$, we have $\co_q(V_\mu)=0$. 
When $-\tilde\mu\in\cP_+$,  we set $$\nu=(-\tilde\mu)^\dagger.$$ 
Then, we may conclude that $\co_q(V_\mu)$ is spanned by
\ba
\zeta_i &=& \sum_{j}
\phi^{(\nu)} (w_j^{(\nu)})\otimes {\tilde t}_{i j}^{(\lambda)},   
\label{irrep}
\na
which are  obviously linearly independent. Furthermore,
\ban x\cdot\zeta_i &=&
(-1)^{[x][\phi^{(\nu)}]}\  \sum_{j} t^{(\nu)}_{j i} (x)\  \zeta_j,
\ \ \ \ \ x\in\Uq.    \nan
Thus $\co_q (V_\mu)\cong W(\nu)$.  More explicitly, 
the isomorphism is given by 
\ba
W(\nu)
\stackrel{(id\otimes S)\delta}{\longrightarrow}
\co_q( W(\nu))
\stackrel{\phi^{(\nu)}\otimes id }{\longrightarrow}
\co_q( V_\mu).  
\na
This completes the proof of the theorem.  
\end{proof}

This result provides an analog of 
the celebrated Borel - Weil theorem for the quantum 
supergroup $\OSPq$.  
For the classical Lie supergroups, the program of developing
a  Bott - Borel - Weil theory was extensively
investigated by Penkov and co - workers\cite{Penkov}, 
although it has not yet been completed so  far as we are aware.
Also,  a quantum Borel-Weil theorem for the covariant and
contravariant tensor representations of quantum $GL(m|n)$
was obtained in \cite{Zhang}.

When $\mu=0$, the theorem implies that 
\ban 
\left\{ f\in\Tq | p\circ f = \epsilon(p) f, \ \ \forall p\in\Up\right\}
& =& \C \epsilon.
\nan
Combining this result with  with Proposition \ref{projective},
we obtain

\bigskip 
\noindent{\bf Corollary }: 
{\em Let $W$ be any finite dimensional $\Uq$-module.
Then, as $\Uq$-modules,
\ban
\co_q(W) &\cong & \epsilon\otimes W.
\nan}

\bigskip
\noindent {\bf Acknowledgement}:  This work is supported  
by the National Science Council (ROC) Grant NSC-87-2112-M008-002. 
Zhang wishes to thank the Department of Physics and 
the Center for Complex Systems at the National Central University 
for the hospitality extended to him during a visit 
from late 1997 to early 1998, when this work was largely 
completed.

\pagebreak

\end{document}